\documentclass[11pt]{article}
\usepackage{amsfonts}
\usepackage{amsmath}
\usepackage{amssymb}
\usepackage{times}
\usepackage{mhequ}
\usepackage{mhenvs}
\usepackage{Symbols}

\def\theequation{\thesection.\arabic{equation}}
\newcommand{\beq}{\begin{equation}}
\newcommand{\eeq}{\end{equation}}
\newcommand{\bea}{\begin{eqnarray}}
\newcommand{\eea}{\end{eqnarray}}
\newcommand{\beas}{\begin{eqnarray*}}
\newcommand{\eeas}{\end{eqnarray*}}

\newcommand{\E}{\mathbb{E}}

\renewcommand{\P}{\mathbb{P}}
\newcommand{\D}{\mathbb{D}}
\renewcommand{\R}{\mathbb{R}}
\newcommand{\loc}{\mathrm{loc}}
\newcommand{\HS}{{\mathrm{HS}}}

\newcommand{\eqlaw}{\stackrel{\mbox{\tiny law}}{=}}

\title{A version of H\"ormander's theorem for the fractional Brownian motion}

\author{Fabrice Baudoin
\\\small Laboratoire de Statistiques et Probabilit\'es
\\\small Universit\'e Paul Sabatier
\\\small 118 Route de Narbonne, Toulouse, France
\\\small fbaudoin@cict.fr
\and Martin Hairer
\\\small Mathematics Institute
\\\small The University of Warwick
\\\small Coventry CV4 7AL, United Kingdom
\\\small hairer@maths.warwick.co.uk }

\begin{document}
\maketitle 

\begin{abstract}
It is shown that the law of an SDE driven by fractional Brownian motion with Hurst
parameter greater than $1/2$ has a smooth density with respect to Lebesgue measure,
provided that the driving vector fields satisfy H\"ormander's condition. The main
new ingredient of the proof is an extension of Norris' lemma to this situation.
\end{abstract}

\tableofcontents

\section{Introduction}

In the celebrated paper \cite{Hor67}, L. H\"ormander gave in 1967
a sufficient (and necessary  in the analytic case) condition for  the hypoellipticity of second order
differential operators. The original proof of H\"ormander was
rather complicated and has been since considerably simplified in
using the theory of pseudo-differential operators.

In 1976, P. Malliavin, uses the deep connection between the theory
of second order differential operators and It\^o's theory of
stochastic differential equations to point out the probabilistic
counterpart of  H\"ormander's theorem. The problem of the
hypoellipticity is essentially equivalent to the problem of the existence of a
smooth density with respect to the Lebesgue measure for the
solution of a stochastic differential equation. The idea of
Malliavin's proof of H\"ormander's theorem is to show that the
It\^o's map associated with a stochastic differential equation is
differentiable in a weak sense and then to show that, under
H\"ormander's conditions, this derivative is non-degenerate. The
stochastic calculus of variations that has been developed in \cite{Mal76}
precisely in order to obtain a probabilistic proof of H\"ormander's theorem is 
now known as Malliavin calculus and
has since then found numerous applications (see \cite{Nua95}).

In the last few years, there have been numerous attempts to define
a notion of solution for differential equations driven by a
fractional Brownian motion. When the Hurst parameter of the
fractional Brownian motion is greater than $1/2$, existence and
uniqueness of the solution are obtained by Z\"ahle in \cite{Zah01}
or Nualart-R\v a\c scanu in \cite{NR02DED}. Let us note that, as a
consequence of the work of Coutin and Qian \cite{CouQia02SAR}, a
notion of solution can actually be well-defined for $H
>\frac{1}{4}$. The problem of the existence and smoothness of the density with
respect to Lebesgue measure for solutions of stochastic
differential equations that are driven by a fractional Brownian
motion with Hurst parameter greater than $1/2$ is solved in some
special cases. In \cite{NouSi06} the existence and smoothness of
the density has been shown in the one-dimensional case by using
Doss-S\"ussman methods.  In \cite{NS05MCS}, the authors prove the
existence of a density under ellipticity assumptions. Finally, in
\cite{HuNu2006}, always under an ellipticity assumption, the
smoothness of the density is proved.

In the present work we prove a version of H\"ormander's theorem
for solutions of differential equations driven by a fractional
Brownian motion. More precisely, we prove that as soon as the
usual H\"ormander conditions are satisfied, there exists a
smooth density for the law of the solutions. One of the major difficulties is
to prove an equivalent for fractional Brownian motion of Norris' lemma
\cite{Nor86SMC}. This lemma, interesting in itself,
quantifies in which way the integrand of a stochastic integral
with respect to fractional Brownian motion has to be small if the
stochastic integral is small.

This article is organised as follows.
In a first section we recall the basics of
Malliavin calculus in the context of fractional Brownian motion, which is the main tool
used in this article. We then prove a version of Norris'
lemma for stochastic integrals with respect to fractional
Brownian motion. This lemma is then used to prove the main result
of this work, that is the existence and the smoothness of the
density under H\"ormander's conditions. Finally, in a last
section, we apply this result to the analysis of the behaviour in small times 
of this density on the diagonal.

{\small

\subsection*{Acknowledgements}

The authors would like to thank L.~Coutin and P.~Friz for stimulating discussions.
}

\section{Malliavin calculus with respect to fractional Brownian
motion} 

Let us first recall some basic facts about Malliavin
calculus with respect to the fractional Brownian motion (for
further details, we refer for instance to \cite{HuNu2006} or
\cite{NS05MCS}).

We consider the Wiener space of continuous paths:
\[
\mathbb{W}^{\otimes d}=\left( \mathcal{C} ( [0,1] , \mathbb{R}^d
), (\mathcal{B}_t)_{0 \leq t \leq 1}, \mathbb{P} \right)
\]
where:
\begin{enumerate}
\item  $\mathcal{C} ( [0,1] , \mathbb{R}^d
)$ is the space of continuous functions $ [0,1] \rightarrow
\mathbb{R}^d$;
\item  $\left( \beta_{t}\right) _{t\geq 0}$ is the coordinate process defined by
$\beta_{t}(f)=f\left( t\right) $, $f \in \mathcal{C} ( [0,1] ,
\mathbb{R}^d )$;
\item  $\mathbb{P}$ is the Wiener measure;
\item  $(\mathcal{B}_t)_{0 \leq t \leq 1}$ is the ($\mathbb{P}$-completed) natural filtration
of $\left( \beta_{t}\right)_{0 \leq t \leq 1}$.
\end{enumerate}
A $d$-dimensional fractional Brownian motion with Hurst parameter
$H \in (0,1)$ is a Gaussian process
\[
B_t = (B_t^1,\ldots,B_t^d), \text{ } t \geq 0,
\]
where $B^1,\ldots,B^d$ are $d$ independent centred Gaussian
processes with covariance function
\[
R\left( t,s\right) =\frac{1}{2}\left(
s^{2H}+t^{2H}-|t-s|^{2H}\right).
\]
It can be shown that such a process admits a continuous version
whose paths are H\"older $p$ continuous, $p<H$.
The fractional Brownian motion with Hurst parameter
$H>\frac{1}{2}$ can be constructed on the Wiener space by a
Volterra type representation. Namely, the process
\begin{equation}\label{Volterra:representation}
B_t = \int_0^t \mathbf{K}(t,s) d\beta_s, t \geq 0
\end{equation}
is a fractional Brownian motion with Hurst parameter $H$, where
\begin{equation*}
\mathbf{K} (t, s)=c_H s^{\frac{1}{2} - H} \int_s^t
(u-s)^{H-\frac{3}{2}} u^{H-\frac{1}{2}} du \;,\qquad t>s.
\end{equation*}
and $c_H$ is a suitable constant.
Let $\mathcal{E}$ be the space of $\mathbb{R}^d$-valued step
functions on $[0,1]$. We denote by $\mathcal{H}$ the closure of
$\mathcal{E}$ for the scalar product:
\[
\langle (\mathbf{1}_{[0,t_1]} , \cdots ,
\mathbf{1}_{[0,t_d]}),(\mathbf{1}_{[0,s_1]} , \cdots ,
\mathbf{1}_{[0,s_d]}) \rangle_{\mathcal{H}}=\sum_{i=1}^d
R(t_i,s_i).
\]
For $\phi,\psi \in \mathcal{H}$, we have
\[
\langle \phi , \psi \rangle_{\mathcal{H}}=H(2H-1)\int_0^1 \int_0^1
\mid s-t \mid^{2H-2} \langle \phi (s) , \psi(s)
\rangle_{\mathbb{R}^d} ds dt.
\]
It can be shown that $L^{1/H} ([0,1], \mathbb{R}^d)
\subset \mathcal{H}$ but that $\mathcal{H}$ also contains
distributions (see for example \cite{PT00IQR}). 

For our purposes, the following representation of the $\CH$-scalar product is useful. 
Let $\CI^\alpha \phi$ denote the fractional
integral of order $\alpha$ of $ \phi$, defined by
\begin{equ}[e:Ifrac]
\CI^{\alpha} \phi(t) ={1 \over \Gamma(\alpha)}\int_0^t (t-s)^{\alpha-1}\phi(s)\,ds \;.
\end{equ}
If one extends $\phi$ and $\psi$ to $\R_+$ by setting $\phi(t) = \psi(t) = 0$
for $t \ge 1$, a straightforward application of Fubini's theorem
shows that one has the identity
\begin{equ}[e:reprH]
\scal{\phi,\psi}_\CH = \scal{\CI^{H-1/2}\phi,\CI^{H-1/2}\psi}_{L^2}\;.
\end{equ}

A $\mathcal{B}_{1}$-measurable real
valued random variable $F$ is said to be cylindrical if it can be
written as
\begin{equation*}
F=f \Bigl( \int_0^{1} \langle h^1_s, dB_s \rangle ,\ldots,\int_0^{1}
\langle h^n_s, dB_s \rangle \Bigr)\;,
\end{equation*}
where $h^i \in \mathcal{H}$ and $f:\mathbb{R}^n \rightarrow
\mathbb{R}$ is a $C^{\infty}$ bounded function. The set of
cylindrical random variables is denoted $\mathcal{S}$. The
Malliavin derivative of $F \in \mathcal{S}$ is the $\mathbb{R}^d$ valued
stochastic process $(\mathbf{D}_t F )_{0 \leq t \leq 1}$ given by
\[
\mathbf{D}_t F=\sum_{i=1}^{n} h^i (t) \frac{\partial f}{\partial
x_i} \left( \int_0^{1} \langle h^1_s, dB_s \rangle,\ldots,\int_0^{1}
\langle h^n_s, dB_s \rangle  \right).
\]
More generally, we can introduce iterated derivatives. If $F \in
\mathcal{S}$, we set
\[
\mathbf{D}^k_{t_1,\ldots,t_k} F = \mathbf{D}_{t_1}
\ldots\mathbf{D}_{t_k} F.
\]
For any $p \geq 1$, the operator $\mathbf{D}^k$ is closable from
$\mathcal{S}$ into $L^p \left( \mathcal{C} ( [0,1] ,
\mathbb{R}^d ) , \mathcal{H}^{\otimes k} \right)$. We denote by
$\mathbb{D}^{k,p}(\mathcal{H})$ the closure of the class of
cylindrical random variables with respect to the norm
\[
\left\| F\right\| _{k,p}=\left( \mathbb{E}\left( F^{p}\right)
+\sum_{j=1}^k \mathbb{E}\left( \left\| \mathbf{D}^j F\right\|
_{\mathcal{H}^{\otimes j}}^{p}\right) \right) ^{\frac{1}{p}},
\]
and
\[
\mathbb{D}^{\infty}(\mathcal{H})=\bigcap_{p \geq 1} \bigcap_{k
\geq 1} \mathbb{D}^{k,p}(\mathcal{H}).
\]
We also introduce the localised spaces $\D^{k,p}_\loc
(\mathcal{H})$ by saying that a random variable $F$ belongs to
$\D^{k,p}_\loc (\mathcal{H})$ if there exists a sequence of sets
$\Omega_n \subset \CB_1$ and random variables $F_n \in \D^{k,p}
(\mathcal{H})$ such that $\Omega_n \uparrow \CC([0,1],\R^d)$
almost surely and such that $F = F_n$ on $\Omega_n$.

We then have the following key result which stems from Theorem
2.1.2 and Corollary 2.1.2. in \cite{Nua95}:
\begin{theorem}\label{theo:dens}
Let $F=(F_1,\ldots,F_n)$ be a $\mathcal{B}_1$-measurable random
vector such that:
\begin{enumerate}
\item For every $i=1,\ldots,n$, $F_i \in \D^{1,2}_\loc(\mathcal{H})$;
\item The matrix
$
\Gamma= \left(  \langle \mathbf{D} F^i , \mathbf{D} F^j
\rangle_{\mathcal{H}}  \right)_{1 \leq i,j \leq n}
$
is invertible almost surely.
\end{enumerate}
Then the law of $F$ has a density with respect to the Lebesgue
measure on $\R^n$. If moreover $F \in \D^\infty (\mathcal{H})$ and, for
every $p >1$,
\[
\mathbb{E} \left( \frac{1}{\mid \det \Gamma \mid ^p} \right) < +
\infty,
\]
then this density is smooth.
\end{theorem}
\begin{remark}
The matrix $\Gamma$ is called the Malliavin matrix\index{Malliavin
matrix} of the random vector $F$.
\end{remark}

\section{Norris' lemma for integrals with respect to fractional Brownian motion}

A main ingredient in many probabilistic proofs of H\"ormander's theorem
(see for example \cite{Nua95}) is Norris' lemma \cite{Nor86SMC}. Loosely speaking, it is a more quantitative
version of the uniqueness property of the semimartingale decomposition, stating
that if a semimartingale is small, then both its bounded variation part and its martingale
part must be small. In other words, the martingale part and the bounded variation
part cannot compensate each other. This section is devoted to the proof of 
\prop{prop:main}, which is a version of Norris' lemma formulated in a framework suitable
for the purposes of this article.

We start by stating a minor variant of a well-known concentration
result for Gaussian measures.
\begin{lemma}\label{lem:conc}
Let $\mu$ be a Gaussian measure on a separable Hilbert space
$\mathcal{H}$ with covariance operator $\Gamma$ and write $T =
\sqrt{\mathbf{tr} \Gamma}$ and $\lambda = \sqrt{\|\Gamma\|}$.
Then, there exists a constant $C$ independent of $\Gamma$ such
that one has the bound
\begin{align*}
\mu \bigl(\bigl|\|x\| - T\bigr| \ge h \bigr) \le C \exp
\Bigl(-{h^2\over 4\lambda^2}\Bigr)\;,
\end{align*}
for every $h \ge 0$.
\end{lemma}
\begin{proof}
We can assume $\mathcal{H}$ to be finite-dimensional; the general
case follows from a simple approximation argument since none of
the constants depends on the dimension of $\mathcal{H}$. Define
$\tilde T = \int \|x\|\,\mu(dx)$. It follows from the
isoperimetric inequality for Gaussian measures
\cite{CirIbrSud76NGS,Tal95COM} that
\begin{align} \label{e:boundOri}
\mu \bigl(\bigl|\|x\| - \tilde T\bigr| \ge h \bigr) \le C \exp
\Bigl(-{h^2\over 2\lambda^2}\Bigr)\;,
\end{align}
for every $h \in \mathbb{R}$. Furthermore, it follows from
\cite[Thm~1.7.1]{Bog98GM} that
\begin{align*}
T^2 - \tilde T^2 \le {\pi^2 \over 4} \lambda^2\;,\quad\text{and
thus}\quad T - \tilde T \le {\pi^2 \lambda^2 \over 4T} \le {\pi^2
\lambda \over 4}\;.
\end{align*}
The claim follows at once.
\end{proof}
We now define the class of Gaussian processes that are of interest
to us. For $H \in (\textstyle{1\over 2},1]$, we say that a Gaussian process $B$ is of
type $H$ if it is centred and the function $f$ defined by
\begin{align*}
f(s,t) = \mathbb{E} \left(B(t)-B(s)\right)^2
\end{align*}
is $\CC^1$ and satisfies
\begin{align}\label{e:boundf}
c_1 |t-s|^{2H} \le f(s,t) \le c_2 |t-s|^{2H} ,\quad |\partial_s
\partial_t f(s,t)| \le c_3 |t-s|^{2H-2} ,
\end{align}
for every pair of times $s,t \in (0,1)$ with $s \neq t$. The following is a direct
consequence of \eref{e:boundf} combined with \ref{lem:conc}.

\begin{lemma}\label{lem:concappl}
Let $B^i$ be i.i.d.\ Gaussian processes of type $H > {1\over 2}$
and let $\delta, N > 0$ be such that $\delta N < 1$. Define the
$\mathbb{R}^N$-valued random variables $X^i$ and the number $T$ by
\begin{align*}
X_n^i = B^i(n\delta) - B^i((n-1)\delta)\;,\quad T^2 = \sum_{n=1}^N
f(n\delta, (n-1)\delta)\;.
\end{align*}
Then, there exist constants $C_1$, $C_2$ such that one has the
bound
\begin{align}
\mathbb{P} \left(\left| |X^i| - T\right| \ge h\right) &\le C_1
\exp \left(- C_2{N h^2 \over (\delta N)^{2H}}\right),
\label{e:bounddiff}\\
\mathbb{P} \left(\left|\langle X^i,X^j \rangle \right| \ge
h^2\right) &\le C_1 \exp \left(- C_2{N h^2 \over (\delta
N)^{2H}}\right), \label{e:boundcov}
\end{align}
for every $h \ge 0$ and every pair $(i,j)$.
\end{lemma}

\begin{proof}
Denote by $\Gamma$ the covariance of $X^i$ (this is independent of
$i$). Then, one has
\begin{align*}
\Gamma_{mn} = {1 \over 2} \int_{m\delta}^{(m+1)\delta}
\int_{n\delta}^{(n+1)\delta} \partial_s \partial_t
f(s,t)\,ds\,dt\;,
\end{align*}
so that $|\Gamma_{mn}| \le C \delta^{2H} (1+|m-n|)^{2H-2}$. This
implies that
\begin{align}\label{e:boundGamma}
\|\Gamma\|_\HS^2 = \sum_{m,n=1}^N |\Gamma_{mn}|^2 \le C
N\delta^{4H} \sum_{k=1}^{N} |k|^{4H-4} \le C \delta^{4H} N^{4H-2}.
\end{align}
Since this is a bound on $\|\Gamma\|^2$, the first inequality
follows from Lemma \ref{lem:conc}. Fix now an arbitrary pair of
indices $i \neq j$. Note that conditional on the value of $X^j$,
the random variable $\langle X^i,X^j \rangle$ is normal with
variance $\langle X^j, \Gamma X^j \rangle$. This motivates the
introduction of the random vector $\bar X^j = \Gamma^{1/2} X^j$
which is Gaussian with covariance $\Gamma^2$. Note also that $\E
\|\bar X^j\|^2 = \mathbf{tr} (\Gamma^2) =
\|\Gamma\|_\HS^2$. We thus have, for any $v > 0$, the
bound
\begin{align*}
\mathbb{P} \left(\left|\langle X^i,X^j \rangle \right| \ge
h^2\right) &= \mathbb{E} \left( \mathbb{P} \left(\left|\langle
X^i,X^j \rangle \right| \ge h^2\,|\,X^j\right)\right)
\le \mathbb{E} \left(  C \exp(-h^4/\|\bar X^j\|^2) \right) \\
&\le C \exp \left(-{h^4\over 4(\|\Gamma\|_\HS +
v)^2}\right)
+ \mathbb{P} \left(|\|\bar X^j\|-\|\Gamma\|_\HS| \ge v\right) \\
&\le C \exp \left(-{h^4\over 8(\|\Gamma\|_\HS^2 +
v^2)}\right)
+ C \exp \left(-{v^2 \over 4 \|\Gamma^2\|}\right) \\
&\le C \exp \left(-{h^4\over 8(\|\Gamma\|_\HS^2 +
v^2)}\right) + C \exp \left(-{v^2 \over 4
\|\Gamma\|_\HS^2}\right)\;.
\end{align*}
Note now that the second inequality is non-trivial only for $h^2
\ge \|\Gamma\|_\HS$, so that we assume that we are in
this situation from now on. Choosing $v^2 = h^2
\|\Gamma\|_\HS$, we get
\begin{align}
\mathbb{P} \left(\left|\langle X^i,X^j \rangle \right| \ge
h^2\right) \le C \exp \left(-h^2 / (16
\|\Gamma\|_\HS)\right)
\end{align}
which, together with \eref{e:boundGamma}, implies the required
bound.
\end{proof}

\begin{remark}
The bounds in \lem{lem:concappl} can be interpreted as saying that the
`coarse-grained quadratic variation' of $B$ on a scale $\delta$ and over
a time interval $t$ behaves like $\delta^{2H-1}t$ to within an error of order
$\delta^H t^H$. Note that the relative magnitude of the error to the average 
value always tends to $0$ as $\delta \to 0$, but that this ratio becomes `worse'
as $H \to 1$.
\end{remark}

We now have the main tools in place to prove the following version of
Norris' lemma. Note that here and in the sequel, we denote by $\|\cdot\|_\alpha$
the $\alpha$-H\"older norm of a function.

\begin{proposition}\label{prop:main}
Let $H \in ({1\over 2},1)$ and let $a$ and $b$ be processes taking
values in $\mathbb{R}$ and $\mathbb{R}^m$ respectively such that
$\mathbb{E} \left(\|a\|_{\tilde H} + \sum_i \|b^i\|_{{\tilde
H}}\right)^p < \infty$ for every $p \ge 1$ and every $\tilde H \in
({1\over 2},H)$. Let
\begin{align*}
y_t = \int_0^t a(s)\,ds +  \int_0^t \langle b(s),dB(s) \rangle \;,
\end{align*}
where the $B_i$ are $m$ i.i.d.\ Gaussian process of type $H$. Then
there exists $q > 0$ such that, for every $p>0$, the estimate
\begin{align}\label{e:mainbound}
\mathbb{P} \left(\|y\|_{L^\infty} < \epsilon \quad\text{and}\quad
\|a\|_{L^\infty} +\|b\|_{L^\infty} > \epsilon^q \right) < C_p
\epsilon^p
\end{align}
holds. The constant $C_p$ depends on $a$, $b$, and $p$ but not on
$\epsilon$.
\end{proposition}

\begin{remark}\label{rem:QV}
Note that we do not require the $B_i$ to be independent of the
processes $a$ and $b$. We also do not require any adaptedness at this stage.
The reason why we will require adaptedness later on is that
equation \eref{e:exprMD} for the Malliavin derivative of the solution does not hold otherwise.
\end{remark}

\begin{remark}
The bound \eref{e:mainbound} actually implies the bound
\begin{equ}[e:main]
\mathbb{P} \left(\|y\|_{L^\infty} < \epsilon\right) \le C_p \eps^p
+ \min \bigl\{ \P(\|a\|_{L^\infty} < \eps^q), \P(\|b\|_{L^\infty}
< \eps^q)\bigr\}\;,
\end{equ}
which will be used repeatedly in the sequel.
\end{remark}

\begin{proof}[of \prop{prop:main}]
The proof consists of two parts. In the first part, we show that
one has a bound of the type
\begin{align}\label{e:boundb}
\mathbb{P} \left(\|y\|_{L^\infty} < \epsilon\quad\text{and}\quad
\|b\|_{L^\infty} > \epsilon^q \right) < C_p \epsilon^p\;.
\end{align}
In the second part, we use this information to show that one has
also
\begin{align}\label{e:bounda}
\mathbb{P} \left(\|y\|_{L^\infty} < \epsilon\quad\text{and}\quad
\|a\|_{L^\infty} > \epsilon^q \right) < C_p \epsilon^p\;.
\end{align}
Combining both bounds then yields \eref{e:mainbound}. 

The key idea to the proof of \eref{e:boundb}, which should be translated to `if $y$ is small
then $b$ must also be small' is the following. Choose two small length scales
$\delta \ll \Delta \ll 1$. Since we assume some regularity on $b$, it is easy to control
the error made by assuming that $b$ is constant on intervals of length $\Delta$.
One then considers the square root of the coarse-grained quadratic variation of $y$ on a scale
$\delta$ over an interval of size $\Delta$ around $t$. This is of course bounded by
$\|y\|_{L^\infty} \delta^{-1/2} \Delta^{1/2}$.
By \rem{rem:QV}, the contribution of the
term including $b$ to this expression is approximately
equal to $\delta^{H-1/2} \Delta^{1/2} |b(t)|$. The contribution of the
term including $a$ on the other hand is bounded by $\|a\|_{L^\infty} \delta^{1/2}\Delta^{1/2}$. 
Summing over all intervals of size $\Delta$ yields  a bound of the type
\begin{equ}
\|b\|_{L^1} \lesssim \delta^{-H} \|y\|_{L^\infty} + \delta^{1-H} \|a\|_{L^\infty}\;,
\end{equ}
from which it is then straightforward to deduce \eref{e:boundb} by making use of the \textit{a priori} bounds
on the H\"older norms of $a$ and $b$.

This argument is of course extremely sloppy, since we have not justified in any way
some of the approximations made and we have not addressed the fact that $b$ takes
values in $\R^m$.
Fix some
small value of $\epsilon > 0$ and fix two small numbers $\delta$
and $\Delta$ such that $1/\delta$ and $1/\Delta$ are integers,
$1/\Delta$ divides $1/\delta$, and such that $\delta \ll \Delta
\ll 1$. We also fix $\tilde H \in ({1\over 2},H)$ to be determined
later. We define $\bar b(t)$ as the stepfunction with steps of
length $\Delta$ approximating $b(t)$, $\bar b(t) =
b(\Delta[t/\Delta])$ and we write $\beta (t) = b(t) - \bar b(t)$.
The stochastic integral of $\beta$ is bounded as follows.
\begin{lemma}\label{lem:bounddb}
There exists a constant $C$ such that, for every $t \ge0$ and
every $s \in [0,\Delta]$, one has
\begin{align*}
\left|\int_t^{t+s} \langle \beta(r),dB(r) \rangle \right| \le C
\|b\|_{\tilde H}\|B\|_{\tilde H}  \Delta^{\tilde H} s^{\tilde
H}\;.
\end{align*}
\end{lemma}

\begin{proof}
Since $s \le \Delta$,  we can assume without loss of generality that $\bar b$ is constant
on the interval $[t,t+s]$. It follows from \cite{You36AIH} that,
for every $\alpha > 1$, there exists a constant $C$ such that
\begin{align*}
\left|\int_t^{t+s} \langle \beta(r),dB(r) \rangle - \langle
\beta(t),B(t+s) - B(t) \rangle \right| \le C \|b\|_{\tilde
H}\|B\|_{\tilde H} s^{2 \tilde H}\;.
\end{align*}
Furthermore, one has $\left|\langle \beta(t),B(t+s) - B(t) \rangle
\right| \le \|b\|_{\tilde H}\|B\|_{\tilde H}\Delta^{\tilde H}
s^{\tilde H}$, so that the result follows at once.
\end{proof}

Denote by $r$ the (integer) ratio $\Delta/\delta$, set $t_n =
\delta n$ and define for $N = 1,\ldots,\Delta^{-1}$ and for $i,j =
1,\ldots,m$ the random variable
\begin{align*}
X_N^{ij} = \sum_{n = (N-1)r}^{Nr-1} \big(B_i(t_{n+1}) -
B_i(t_n)\bigr)\big(B_j(t_{n+1}) - B_j(t_n)\bigr)\;.
\end{align*}
Note that $T^2 := \E X_N^{ii} = r f(\delta) \approx \Delta
\delta^{2H-1}$. With these notations and using Lemma
\ref{lem:bounddb}, we have, for $n \in [(N-1)r, Nr-1]$, the
relation
\begin{align*}
|\langle b(N\Delta), B(t_{n+1}) - B(t_n) \rangle | &\le |y_{t_{n+1}} - y_{t_{n}}| + \|a\|_{L^\infty} \delta + \left|\int_{t_n}^{t_{n+1}} \langle \beta(s),dB(s) \rangle \right| \\
&\le 2\|y\|_{L^\infty} + \|a\|_{L^\infty} \delta + C
\|b\|_{\tilde H}\|B\|_{\tilde H}  \Delta^{\tilde H} \delta^{\tilde
H}\;.
\end{align*}
Taking squares on both sides, summing from $n = (N-1)r$ to $n =
Nr-1$ and taking square roots, we get
\begin{align}\label{e:boundDelta}
\sqrt{\sum_{i,j} b_i(N\Delta)b_j(N\Delta) X_N^{ij}} \le
\Delta^{1/2}\delta^{-1/2}\bigl(2\|y\|_{L^\infty} +
\|a\|_{L^\infty} \delta + C \|b\|_{\tilde H}\|B\|_{\tilde H}
\Delta^{\tilde H} \delta^{\tilde H}\bigr)\;.
\end{align}
At this point, it is convenient to introduce quantities $Y_N^i =
(X_N^{ii})^{1/2}$ and $Y_N^{ij} = |X_N^{ij}|^{1/2}$. With this
notation, we get
\begin{align*}
\sum_{i=1}^m |b_i(N\Delta)| Y_N^{i} &\le C \sum_{i \neq j}
Y_N^{ij} \sqrt{|b_i(N\Delta)b_j(N\Delta)|} +
C\Delta^{1/2}\delta^{-1/2}\|y\|_{L^\infty} \\
& \quad + C \Delta^{1/2} \delta^{1/2} \|a\|_{L^\infty}  + C
\Delta^{\tilde H + 1/2} \delta^{\tilde H-1/2} \|b\|_{\tilde
H}\|B\|_{\tilde H}\;.
\end{align*}
Summing over $N$ yields
\begin{align}\label{e:boundsumb}
\sum_{i=1}^m \sum_{N=1}^{\Delta^{-1}} |b_i(N\Delta)| Y_N^{i} &\le
C \sum_{i \neq j} \sum_{N=1}^{\Delta^{-1}} Y_N^{ij}
\sqrt{|b_i(N\Delta)b_j(N\Delta)|} +
C\Delta^{-1/2}\delta^{-1/2}\|y\|_{L^\infty} \\
& \quad + C \Delta^{-1/2} \delta^{1/2} \|a\|_{L^\infty}  + C
\Delta^{\tilde H - 1/2} \delta^{\tilde H-1/2} \|b\|_{\tilde
H}\|B\|_{\tilde H}\;.
\end{align}
Note now that by Lemma \ref{lem:conc} one can hope that $Y_N^{i}
\approx T$ for every $i$ and $N$. Therefore, the left hand side is
very close to $T/\Delta$ times the sum of the $L^1$ norms of the
$b_i$. More precisely, one has for every $i$, the bound
\begin{align*}
\Bigl|\Delta \sum_{N=1}^{\Delta^{-1}} |b_i(N\Delta)| -
\|b_i\|_{L^1}\Bigr| \le \|b_i\|_{\tilde H} \Delta^{\tilde H}\;,
\end{align*}
so that, multiplying \eref{e:boundsumb} by $\Delta/T$, we get the
bound
\begin{equs}
\sum_{i=1}^m \|b_i\|_{L^1} &\le  C\Delta^{\tilde H} \|b\|_{\tilde
H}  + C\delta^{-H}\|y\|_{L^\infty}
+ C\Delta^{\tilde H}\delta^{\tilde H - H}\|b\|_{\tilde H}\|B\|_{\tilde H}\\
&\quad + C\delta^{1-H}\|a\|_{L^\infty} + C\Delta^{1/2} \delta^{1/2-H} \sum_{i=1}^m \sum_{N=1}^{\Delta^{-1}} |b_i(N\Delta)| |Y_N^i-T|\\
&\quad + C \Delta^{1/2}\delta^{1/2-H} \sum_{i \neq j} \sum_{N=1}^{\Delta^{-1}} Y_N^{ij} \sqrt{|b_i(N\Delta)b_j(N\Delta)|} \\
&\le C\Delta^{\tilde H} \|b\|_{\tilde H}  +
C\delta^{-H}\|y\|_{L^\infty}
+ C\Delta^{\tilde H}\delta^{\tilde H - H}\|b\|_{\tilde H}\|B\|_{\tilde H}\\
&\quad + C\delta^{1-H}\|a\|_{L^\infty} + C\Delta^{1/2}
\delta^{1/2-H} \|b\|_{\tilde H} \sum_{i,j=1}^m
\sum_{N=1}^{\Delta^{-1}} |Y_N^{ij}-T \delta_{ij}|\;,
\end{equs}
where we used $\delta_{ij}$ to denote the Kronecker delta. At this
point, we note that, for every $\gamma \le 1$, one has the
interpolation inequality
\begin{align}\label{e:interp}
\|b\|_{L^\infty} \le C \bigl(\gamma \|b\|_{\tilde H} +
\gamma^{-1/\tilde H} \|b\|_{L^1}\bigr)\;.
\end{align}
Therefore
\begin{align*}
\|b\|_{L^\infty} &\le C \gamma^{-1/\tilde H}\Delta^{\tilde H} \|b\|_{\tilde H}  + C\gamma^{-1/\tilde H} \delta^{-H}\|y\|_{L^\infty} + C\gamma^{-1/\tilde H} \delta^{1-H}\|a\|_{L^\infty} \\
&\quad + C\gamma^{-1/\tilde H} \Delta^{\tilde H}\delta^{\tilde H -
H}\|b\|_{\tilde H}\|B\|_{\tilde H}
+ C \gamma \|b\|_{\tilde H}\\
&\quad+ C\gamma^{-1/\tilde H} \Delta^{1/2} \delta^{1/2-H}
\|b\|_{\tilde H} \sum_{i,j=1}^m \sum_{N=1}^{\Delta^{-1}}
|Y_N^{ij}-T\delta_{ij}| \;.
\end{align*}
We now make the following choices for $\gamma$, $\delta$, and
$\Delta$:
\begin{align*}
\gamma \approx \epsilon^{{H(1-H) \over (1+H)(2-H)}}\;,\quad \delta
\approx \epsilon^{1 \over H(2-H)}\;,\quad \Delta\approx
\epsilon^{1-H \over H(2-H)}\;.
\end{align*}
Note that, provided $\epsilon$ is small, one has indeed $\delta
\ll \Delta \ll 1$. With these choices, and by choosing $\tilde H$
sufficiently close to $H$, we see that there exists a constant
$\alpha > 0$ and a constant $C > 0$ such that
\begin{align*}
\|b\|_{L^\infty} &\le C\epsilon^{\alpha} \left(\|b\|_{\tilde H}
\left(1 + \|B\|_{\tilde H} +  \sum_{i,j,N}
{|Y_N^{ij}-T\delta_{ij}|\over \delta^{1/2} \Delta^{H-3/2}}
\right) + \epsilon^{-1}\|y\|_{L^\infty} + \|a\|_{L^\infty}\right)\\
&\le C\epsilon^{\alpha} \left(1 + \|b\|_{\tilde H}^2 +
\|B\|_{\tilde H}^2 +
\left(\sum_{i,j,N}{|Y_N^{ij}-T\delta_{ij}|\over \delta^{1/2}
\Delta^{H-3/2}}\right)^2 + {\|y\|_{L^\infty}\over \epsilon} +
\|a\|_{L^\infty}\right)\;.
\end{align*}
Actually, the constant $\alpha$ can be brought arbitrarily close
to ${H(1-H)\over(1+H)(2-H)} \ge {2\over 9}(1-H)$. Note that Lemma
\ref{lem:concappl} yields the bound
\begin{align*}
\mathbb{P} \left(|Y_N^{ij}-T\delta_{ij}| > h\right) \le C \exp
\left(-c {h^2 \over \delta \Delta^{2H-1}}\right)\;,
\end{align*}
for every possible value of $i$, $j$, and $N$. This immediately
implies
\begin{align}\label{e:boundvar}
\mathbb{P} \left(\sum_{i,j,N} {|Y_N^{ij}-T\delta_{ij}|\over
\delta^{1/2} \Delta^{H-3/2}}> h\right) \le {C\over \Delta} \exp
\left(-c h^2\right)\;.
\end{align}
Therefore, there exists a constant $c$ such that, for $\epsilon$
small enough, one has the bound
\begin{equs}\label{e:bound1}
\mathbb{P} \left(\|b\|_{L^\infty} > \epsilon^{\alpha/2} \quad
\text{and}\quad\|y\|_{L^\infty} < \epsilon \right) & \le
\mathbb{P} \left(\sum_{i,j,N} {|Y_N^{ij}-T\delta_{ij}|\over
\delta^{1/2} \Delta^{H-3/2}}>
c\epsilon^{-\alpha/4}\right) \\
  & + \mathbb{P} \left(\|b\|_{\tilde H}^2 + \|B\|_{\tilde H}^2 +
\|a\|_{L^\infty} > c \epsilon^{-\alpha/2} \right) \le C_p
\epsilon^p\;.
\end{equs}
The last inequality is obtained by combining \eref{e:boundvar}
with the \textit{a priori} bounds on the processes $a$, $b$, and
$B$. We now turn to the proof of \eref{e:bounda}. Fix again some
small value of $\Delta$ to be determined later. We then have for
every $t \in [0,1-\Delta]$ the inequality
\begin{align*}
\left|\int_t^{t+\Delta}a(s)\,ds\right| &\le 2\|y\|_{L^\infty} +
\left|\langle
b(t), B(t+\Delta)-B(t) \rangle \right| + C \|b\|_{\tilde H} \|B\|_{\tilde H} \Delta^{2\tilde H}\\
&\le 2\|y\|_{L^\infty} + \|b\|_{L^\infty} \|B\|_{\tilde
H}\Delta^{\tilde H} + C \|b\|_{\tilde H} \|B\|_{\tilde H}
\Delta^{2\tilde H}\;.
\end{align*}
It is easy to show that, similarly to \eref{e:interp} one has the
inequality
\begin{align*}
\|a\|_{L^\infty} \le 2\Delta^{-1}\sup_{t \in
[0,1-\Delta]}\Bigl|\int_t^{t+\Delta}a(s)\,ds\Bigr| +
2\Delta^{\tilde H}\|a\|_{\tilde H}\;,
\end{align*}
so that
\begin{align*}
\|a\|_{L^\infty} &\le 2\Delta^{-1} \|y\|_{L^\infty} + C \|b\|_{L^\infty} \|B\|_{\tilde H}\Delta^{\tilde H-1} + C \|b\|_{\tilde H} \|B\|_{\tilde H} \Delta^{2\tilde H-1} + 2\Delta^{\tilde H}\|a\|_{\tilde H}\\
&\le C \bigl(\Delta^{-1}  \bigl(\|y\|_{L^\infty} +
\|b\|_{L^\infty}^2\bigr) + \Delta^{2\tilde H-1}
\bigl(\|B\|_{\tilde H}^2 + \|b\|_{\tilde H}^2 + \|a\|_{\tilde
H}\bigr)\bigr)\;.
\end{align*}
At this point, we choose $\Delta \approx \epsilon^{\alpha / H}$
with $\alpha$ as in \eref{e:bound1}. This implies (by choosing as
before $\tilde H$ sufficiently close to $H$) that there exists
$\beta > 0$ such that
\begin{align*}
\|a\|_{L^\infty} \le C \epsilon^{\beta}\left({\|y\|_{L^\infty}
\over \epsilon} + {\|b\|_{L^\infty}^2 \over \epsilon^\alpha} +
\|B\|_{\tilde H}^2 + \|b\|_{\tilde H}^2 + \|a\|_{\tilde
H}\right)\;.
\end{align*}
Therefore, for $\epsilon$ small enough, there exists a constant
$c$ such that
\begin{align*}
\mathbb{P} \left(\|a\|_{L^\infty} > \epsilon^{\beta/2} \quad
\text{and}\quad\|y\|_{L^\infty} < \epsilon\right) & \le \mathbb{P}
\left(\|b\|_{L^\infty} > \epsilon^{\alpha/2}
\quad\text{and}\quad\|y\|_{L^\infty}
< \epsilon \right)\\
&+\mathbb{P} \left(\|B\|_{\tilde H}^2 + \|b\|_{\tilde H}^2 +
\|a\|_{\tilde H} > c\eps^{-\beta/2} \right) \le C_p \epsilon^p\;,
\end{align*}
for arbitrary values of $p$. The last inequality is obtained by
combining \eref{e:bound1} with the \textit{a priori} bounds on the
processes $a$, $b$, and $B$. This concludes the proof of
\eref{e:bounda} and thus of Proposition \ref{prop:main}.
\end{proof}

\section{Existence and smoothness of the density under H\"ormander's
condition}

We now arrive to the heart of our study and are interested in the
study of the existence and regularity for the density of solutions
of stochastic differential equations on $\mathbb{R}^n$
\begin{equation}
\label{SDEmalliavin} X^{x}_t =x +\int_0^t V_0 (X^x_s)ds+
\sum_{i=1}^d \int_0^t V_i (X^{x}_s) dB^i_s
\end{equation}
where the $V_i$'s are $C^{\infty}$-bounded vector fields on
$\mathbb{R}^n$ and $B$ is the $d$ dimensional fractional Brownian
motion defined by (\ref{Volterra:representation}). For this type
of equations,  existence and uniqueness of the solution have been
investigated by many authors (for instance  Nualart-Rasc\^anu in
\cite{NR02DED}). It has also been shown very recently in \cite{HuNu2006}
that the law of the solution to \eref{SDEmalliavin}
possesses a smooth density with respect to the Lebesgue measure in the \textit{elliptic}
case.
Recall first the following \textit{a priori} bound on the
solutions to \eref{SDEmalliavin}. (See for example \cite{NR02DED} for a proof.)

\begin{lemma}\label{lem:apriori}
For every $\gamma < H$ and every $p > 0$,
$\mathbb{E} \left( \| X \|_\gamma^p \right)<+\infty$.
\end{lemma}

Let us now denote by $\Phi$ the stochastic flow associated with
equation (\ref{SDEmalliavin}), that is $\Phi_t (x)=X_t^x$. From
\cite{HuNu2006}, we can deduce:
\begin{lemma} \label{differential SDE}
The map $\Phi_t$ is $C^1$ and the first variation process defined
by
\[
\mathbf{J}_{0 \rightarrow t}=\frac{\partial \Phi_t}{\partial x},
\]
satisfies the following equation:
\[
\mathbf{J}_{0 \rightarrow t}=\mathbf{Id}_{\mathbb{R}^n}+\int_0^t
DV_0(X_s^x)\mathbf{J}_{0 \rightarrow s}ds+ \sum_{i=1}^d \int_0^t
DV_i (X^{x}_s)\mathbf{J}_{0 \rightarrow s} dB^i_s.
\]
and, for every $p>1$,
\[
\mathbb{E} \left( \parallel \mathbf{J}_{0 \rightarrow 1}^{-1}
\parallel ^{p} \right) < + \infty.
\]
Furthermore, for every $i=1,\ldots,n$, $t>0$, and $x \in \R^n$, $X_t^{x,i} \in
\mathbb{D}^{\infty}(\mathcal{H})$. Moreover,
\begin{equ}[e:exprMD]
\mathbf{D}^j_s X_t^{x}= \mathbf{J}_{0 \rightarrow t}
\mathbf{J}_{0 \rightarrow s}^{-1} V_j (X_s) ,~~j=1,\ldots,d, ~~ 0\leq
s \leq t,
\end{equ}
where $\mathbf{D}^j_s X^{x,i}_t $ is the $j$-th component of
$\mathbf{D}_s X^{x,i}_t$.
\end{lemma}

We can now turn to our version of H\"ormander's theorem for
stochastic differential equations that are driven by a fractional
Brownian motion with Hurst parameter $H>1/2$. If $I=(i_1,\ldots,i_k)
\in \{ 0,\ldots, d \}^k$, we denote by $V_I$ the Lie commutator
defined by
\[
V_I = [V_{i_1},[V_{i_2},\ldots,[V_{i_{k-1}}, V_{i_{k}}]\ldots].
\]
We also define the sets of vector fields
\begin{equ}
\CV_n = \bigl\{V_I\,,\,I \in \{0,\ldots,d\} \times
\{1,\ldots,d\}^{n-1}\bigr\}\;, \quad \bar \CV_n = \bigcup_{k=0}^n
\CV_k\;.
\end{equ}
With these notations, the main result of this article is the
following:

\begin{theorem}\label{theo:main}
Assume that, at some $x_0 \in \mathbb{R}^n$, there exists $N$ such
that
\begin{equation}
\label{Hormander assumption} \mathbf{span} \{ V(x_0), \, V\in
\bar\CV_N \} =\mathbb{R}^n\;.
\end{equation}
Then, for any $t>0$, the law of the random variable $X_t^{x_0}$
has a smooth density with respect to the Lebesgue measure on
$\mathbb{R}^n$.
\end{theorem}

Given the results from the previous sections, the proof of this
theorem is by now quite standard and follows closely the argument
given for instance in \cite{Bismut81}, \cite{Nor86SMC} or
\cite{Nua95}. The main difference is that it is not \textit{a priori} obvious
how to relate the $L^\infty$ bounds obtained in \prop{prop:main}
to the fractional Sobolev norms appearing in the statement of \theo{theo:dens}.

\begin{lemma}\label{lem:normH}
Let $H > 1/2$ and let $\CH$ be defined as above. Then, for every
$\gamma > H - 1/2$ there exists a constant $C$ such that
\begin{equ}
\|f\|_\CH \ge C{\|f\|_{L^\infty}^{3 + 1/\gamma} \over
\|f\|_\gamma^{2 + 1/\gamma}}\;,
\end{equ}
for every continuous function $f \in \CH$. Here, $\|f\|_\gamma$ denotes as before
the $\gamma$-H\"older norm of $f$. 
\end{lemma}

\begin{proof}
Let $\CD^\alpha  f$ denote the fractional
derivative of order $\alpha$ of $f$, defined by
\begin{equ}[e:Dfrac]
\CD^{\alpha} f(t) ={1 \over \Gamma(1-\alpha)}{d \over dt} \int_0^t
(t-s)^{-\alpha} f(s)\,ds \;.
\end{equ}
We also introduce the operator $\CD_-^\alpha$ defined by
\begin{equ}
\CD_-^{\alpha} f(t) =-{\alpha \over \Gamma(1-\alpha)}\int_t^\infty
(s-t)^{-\alpha-1}\bigl(f(s)- f(t)\bigr)\,ds \;,
\end{equ}
which is nothing but the adjoint of $\CD^\alpha$ in $L^2(\R_+)$.

Since $\CI^\alpha$ and $\CD^\alpha$ are each other's inverse \cite{SKM93FID}, 
\eref{e:reprH} implies by Cauchy-Schwartz that
\begin{equ}[e:CS]
\bigl|\scal{f,g}_{L^2}\bigr| = \bigl|\scal{\CI^{H-1/2} f,\CD_-^{H-1/2} g}_{L^2}\bigr| \le \|f\|_\CH \|\CD_-^{H-1/2} g\|\;.
\end{equ}
The problem with \eref{e:CS} is that we would like to apply it to
a function $f$ which is $\gamma$-H\"older continuous on $[0,1]$,
but does not necessarily vanish at either $0$ or $1$, so that
$\CD_-^{H-1/2} f$ does in general not belong to $L^2$. If we
define however $h(t) = t^\gamma (1-t)^\gamma$ for $t \in [0,1]$
and $h(t) = 0$ for $t \ge 1$, then $\|fh\|_\gamma \le C
\|f\|_\gamma$, but $fh$ vanishes at $0$ and at $1$. In particular,
this implies that $ \|\CD_-^{H-1/2} fh \| \le C \|f\|_\gamma$, so
that
\begin{equ}
\|f\|_\CH \ge C { \int_0^1 t^\gamma (1-t)^\gamma f^2(t)\,dt \over
\|f\|_\gamma}\;,
\end{equ}
for some constant $C$. On the other hand, it is a straightforward
calculation to check that
\begin{equ}
\int_0^1 t^\gamma (1-t)^\gamma f^2(t)\,dt \ge C{\|f\|_\infty^{3 +
1/\gamma} \over \|f\|_\gamma^{1 + 1/\gamma}}\;,
\end{equ}
which implies the desired result.
\end{proof}

\begin{corollary}\label{cor:boundH}
Let $y$ be a random process with sample paths that are almost
surely $\gamma$-H\"older continuous for some $\gamma > H-1/2$.
Then, there exists an exponent $\alpha$ such that
\begin{equ}
\P \bigl(\|y\|_\CH < \eps \bigr) \le \P \bigl(\|y\|_{L^\infty} <
\eps^{\alpha} \bigr) + \P \bigl(\|y\|_{\gamma} > \eps^{-\alpha}
\bigr)\;,
\end{equ}
for every $\eps$ sufficiently small.
\end{corollary}

\begin{proof}
It follows from \lem{lem:normH} that
\begin{equ}
\P \bigl(\|y\|_{\CH} < \eps\bigr) \le  \P
\bigl(C\|y\|_{L^\infty}^{3 + 1/\gamma} \|y\|_\gamma^{-2-1/\gamma}
< \eps\bigr) \;.
\end{equ}
For an arbitrary pair of positive random variables $X$ and $Y$,
one always has
\begin{equ}
\P (X/Y < \eps) \le \P\bigl(X < \eps^{1-\alpha}\bigr) + \P \bigl(Y
> \eps^{-\alpha}\bigr)\;,
\end{equ}
so that the claim follows by taking $\alpha$ small enough.
\end{proof}

This provides us with the necessary tools to complete the

\begin{proof}[of \theo{theo:main}]
We shall show that $X_1^{x_0}$ admits a smooth density with
respect to the Lebesgue measure, by using the Malliavin covariance
matrix $\Gamma_1$ associated with $X_1^{x_0}$. Note that we can
consider the case $t=1$ without any loss of generality by
rescaling the vector fields $V_i$ appropriately.

Let $\Gamma_1$ be the Malliavin covariance matrix associated with
$X_1^{x_0}$. By definition, we have
\[
\Gamma_1= \bigl(   \langle \mathbf{D} X^{i,x_0}_1 , \mathbf{D}
X^{j,x_0}_1 \rangle_{\mathcal{H}}  \bigr)_{1 \leq i,j \leq n}.
\]

To show that $X_1^{x_0}$ has a a smooth density, it suffices to
show that with probability one $\Gamma_1$ is invertible and that
for every $p>1$,
\[
\mathbb{E} \left( \frac{1}{\mid \det \Gamma_1 \mid ^p} \right) <
+\infty.
\]
From Lemma \ref{differential SDE},
\[
\mathbf{D}^j_s X_1^{x_0}= \mathbf{J}_{0 \rightarrow 1}
\mathbf{J}_{0 \rightarrow s}^{-1} V_j (X_s) ,~~j=1,\ldots,d, ~~ 0\leq
s \leq 1.
\]
Therefore,
\[
\Gamma_1=  H(2H-1)\mathbf{J}_{0 \rightarrow 1} \int_0^1 \int_0^1
\mathbf{J}_{0 \rightarrow u}^{-1} V (X^{x_0}_u)  V (X^{x_0}_v) ^T
\bigl(\mathbf{J}_{0 \rightarrow v}^{-1}\bigr)^T \mid u-v
\mid^{2H-2}\, du\, dv\, \mathbf{J}_{0 \rightarrow 1}^T\;,
\]
where $V$ denotes the $n \times d$ matrix $(V_1 \ldots V_d)$. Since
$\mathbf{J}_{0 \rightarrow 1}$ is almost surely invertible with
inverse in $L^p$, $p>1$, in order to show that $\Gamma_1$ is
invertible with probability one, it is enough to check that with
probability one, the matrix
\[
C_1=\int_0^1 \int_0^1 \mathbf{J}_{0 \rightarrow u}^{-1} V
(X^{x_0}_u)  V (X^{x_0}_v) ^T \bigl(\mathbf{J}_{0 \rightarrow
v}^{-1}\bigr)^T \mid u-v \mid^{2H-2} du dv
\]
is invertible and satisfies for every $p>1$,
\begin{equ}[e:boundCinv]
\mathbb{E} \left( \frac{1}{\mid \det C_1 \mid ^{p}} \right)
<+\infty.
\end{equ}

By using Proposition \ref{prop:main}, the idea is now to control
the smallest eigenvalue of $C_1$ by showing that it can not be too
small. More precisely, recall (Lemma 2.3.1. in \cite{Nua95}) that
if for any $p \ge 2$, there exists $\epsilon_0 (p)$ such that for
every $\epsilon \le \epsilon_0 (p)$,
\[
\sup_{\parallel v \parallel =1} \mathbb{P} \left( \scal{v, C_1 v}
\le \epsilon \right) \le \epsilon^p,
\]
then $C_1$ is invertible with probability one and
(\ref{e:boundCinv}) holds for every $p>1$. We thus want to
estimate $ \mathbb{P} \left( \scal{v, C_1 v} \le \epsilon
\right)$. Let us observe that
\begin{equs}
\scal{v, C_1 v}&=\sum_{j=1}^d \int_0^1 \int_0^1 \mid s-t
\mid^{2H-2} \langle v , (\Phi_s^*V_j)(x_0) \rangle \langle v ,
(\Phi_t^*V_j)(x_0) \rangle \,ds\,dt\\
&= \sum_{j=1}^d \|\scal{v,(\Phi_\cdot^*V_i)(x_0)}\|_\CH^2\;.
\end{equs}
Fix now an (arbitrarily large) value $p > 0$. It follows from
\lem{lem:apriori} and \cor{cor:boundH} that there exists $\alpha >
0$ such that
\begin{equ}[e:startbound]
\mathbb{P} \bigl( \scal{v, C_1 v} \le \eps \bigr) \le C\eps^p +
\min_{i=1,\ldots,d } \mathbb{P} \left(\left\| \left\langle v ,
(\Phi_\cdot^*V_i)(x_0)\right\rangle \right\|_{L^\infty} \le
\epsilon^\alpha \right)\;.
\end{equ}
Note now that if $V$ is an arbitrary bounded vector field with
bounded derivatives, the chain rule reads
\begin{equ}
(\Phi_t^*V)(x_0) = \int_0^t \scal[b]{y, \left( \Phi_s^{\ast} [V_0
, V ] \right) (x_0)}\, ds+\sum_{j=1}^d \int_0^t \scal[b]{y, \left(
\Phi_s^{\ast} [V_j , V ] \right) (x_0) }\, dB^j_s \;.
\end{equ}
It thus follows from Proposition \ref{prop:main} and
\cor{cor:boundH} that there exists $\alpha$ such that
\begin{equ}[e:iterbound]
\P \bigl(\|\scal{v,(\Phi_\cdot^*V)(x_0)}\|_{L^\infty} <
\eps\bigr) \le C\eps^p + \min_{i=0,\ldots,d} \P
\bigl(\|\scal{v,(\Phi_\cdot^*[V_i,V])(x_0)}\|_{L^\infty} <
\eps^\alpha\bigr)
\end{equ}
Consider now the integer $N$ from the assumption. Combining
\eref{e:startbound} with \eref{e:iterbound}, we see that there
exists $\alpha > 0$ such that
\begin{equ}
\mathbb{P} \bigl( \scal{v, C_1 v} \le \eps \bigr) \le C\eps^p +
\min_{V \in \bar \CV_N} \mathbb{P} \left(\left\| \left\langle v ,
(\Phi_\cdot^*V)(x_0)\right\rangle \right\|_{L^\infty} \le
\epsilon^\alpha \right)\;.
\end{equ}
for all $\epsilon$ small enough. On the other hand, we know by
assumption that $\{ V(x_0)\,,\, V \in \bar \CV_N\}$ spans all of
$\R^n$, so that there exists some $V \in \bar \CV_N$ such that
$\scal{v, V(x_0)} \neq 0$. Therefore, one has $\mathbb{P} \bigl(
\scal{v, C_1 v} \le \eps \bigr) \le C\eps^p$ for all $\eps$
sufficiently small, which is the required bound.
\end{proof}

\section{Asymptotics of the density in small times}

In order to obtain asymptotics of the density of hypoelliptic
diffusions on the diagonal in small times, one method consists to
approximate the diffusion by the lift of the Brownian motion in a
nilpotent Lie group that is called a Carnot group (see
\cite{Bau2005} or \cite{BenArous(1989a)}). In a recent work
\cite{BauCou2006}, Baudoin and Coutin have introduced and studied
fractional Brownian motions on Carnot groups. We shall see that
the solution of a stochastic differential equation driven by
fractional Brownian motions can, under H\"ormander's type
assumptions, be approximated by fractional Brownian motions on
Carnot groups. From this approximation, we will deduce an
asymptotic development of the density in small times.

We recall first the notion of Carnot group  (see e.g.
\cite{Bau2005}) and the main results that are obtained in
\cite{BauCou2006}.
\begin{definition}
A Carnot group of step (or depth) $N$ is a simply connected Lie
group $\mathbb{G}$ whose Lie algebra can be decomposed as
\[
\mathcal{V}_{1}\oplus\ldots\oplus \mathcal{V}_{N},
\]
where
\[
\lbrack \mathcal{V}_{i},\mathcal{V}_{j}]=\mathcal{V}_{i+j}
\]
and
\[
\mathcal{V}_{s}=0,\text{ for }s>N.
\]
\end{definition}

Notice that the vector space $\mathcal{V}_{1}$, which is called
the basis of $\mathbb{G}$, Lie generates $\mathfrak{g}$, where
$\mathfrak{g}$ denotes the Lie algebra of $\mathbb{G}$. It is
possible to show that for every $ N \ge 1$, up to an isomorphism,
there exists exactly one $N$-step nilpotent Carnot group with
basis $\mathbb{R}^d$. This group shall be denoted by $\mathbb{G}_N
(\mathbb{R}^d)$.

Since a Carnot group $\mathbb{G}$ is $N$-step nilpotent and simply
connected, the exponential map is a diffeomorphism. On
$\mathfrak{g}$ we can consider the family of linear operators
$\delta_{t}:\mathfrak{g} \rightarrow \mathfrak{g}$, $t \geq 0$
which act by scalar multiplication $t^{i}$ on $\mathcal{V}_{i} $.
These operators are Lie algebra automorphisms due to the grading.
The maps $\delta_t$ induce Lie group automorphisms $\Delta_t
:\mathbb{G} \rightarrow \mathbb{G}$ which are called the canonical
dilations of $\mathbb{G}$. Let us now take a basis $U_1,\ldots,U_d$
of the vector space $\mathcal{V}_1$. The vectors $U_i$ can be
seen as left invariant vector fields on $\mathbb{G}$ so that we
can consider the following stochastic differential equation on
$\mathbb{G}$:
\begin{equation}
\label{SDEcarnot} dX_t =\sum_{i=1}^d \int_0^t U_i (X_s)  dB^i_s,
\text{ } t \geq 0,
\end{equation}
which is easily seen to have a unique solution associated with the
initial condition $X_0=1_{\mathbb{G}}$. The driving process
$(B_t)_{t \ge 0}$ is here a fractional Brownian motion with Hurst
parameter $H>1/2$. The process $(X_t)_{t \ge 0}$ is called the
lift of $(B_t)_{t \ge 0}$ in the group $\mathbb{G}$. For this
equation, the assumptions of Theorem~\ref{theo:main} are
obviously satisfied, so that we have a smooth density for $X_t$,
$t>0$, with respect to the Haar measure of $\mathbb{G}$.

We have then the global scaling property
\begin{equ}
(X_{ct})_{t \geq 0} \eqlaw (\Delta_{c^H} X_{t})_{t \geq 0}.
\end{equ}
(See \cite{BauCou2006} for a proof.)
This scaling property leads directly to the following value at
$1_{\mathbb{G}}$ of the density $\tilde{p}_t$ of $X_{t}$ with
respect to the Haar measure of $\mathbb{G}$:
\[
\tilde{p}_t \left( 1_{\mathbb{G}} \right)= \frac{C}{t^{DH}},
\text{ } t>0,
\]
where $C>0$ and $D=\sum_{i=1}^N i \dim \mathcal{V}_i$. From this,
we will see how to deduce asymptotics in small times for any
\textit{hypoelliptic} stochastic differential equation driven by
fractional Brownian motions.

From now on, we consider $d$ vector fields $V_i: \mathbb{R}^n
\rightarrow \mathbb{R}^n$ which are $C^{\infty}$ bounded  and
shall always assume that the following assumption is satisfied.

\noindent\textbf{Strong H\"ormander's
Condition}\index{H\"ormander's condition}: \textit{For every} $x
\in \mathbb{R}^n$, \textit{we have}:
\[
\mathbf{span} \{ V_I (x), I \in \cup_{k \geq 1} \{1,\ldots,d\}^k \}
=\mathbb{R}^n.
\]

We recall that if $I=(i_1,\ldots,i_k) \in \{ 1,\ldots, d \}^k$ is a
word, we denote by $V_I$ the commutator defined by
\[
V_I = [V_{i_1},[V_{i_2},\ldots,[V_{i_{k-1}}, V_{i_{k}}]\ldots].
\]

Let us introduce some concepts of differential
geometry. The set of linear combinations with smooth coefficients
of the vector fields $V_1,\ldots,V_d$ is called the differential
system\index{Differential system} (or sheaf) generated by these
vector fields. It shall be denoted by $\mathcal{D}$ in the sequel.
Notice that $\mathcal{D}$ is naturally endowed with a
$\mathcal{C}_{\infty} ( \mathbb{R}^n , \mathbb{R})$-module structure. For $x
\in \mathbb{R}^n$, we denote
\[
\mathcal{D} (x)=\{ X(x), X \in \mathcal{D} \}.
\]
If the integer $\dim \mathcal{D} (x)$ does not depend on $x$, then
$\mathcal{D}$ is said to be a distribution. The Lie brackets of
vector fields in $\mathcal{V}$ generate a flag of differential
systems,
\[
\mathcal{D} \equiv \mathcal{D}^1 \subset \mathcal{D}^2 \subset \cdots \subset
\mathcal{D}^k \subset \cdots ,
\]
where $\mathcal{D}^k$ is recursively defined by the formula
\[
\mathcal{D}^k=\mathcal{D}^{k-1} + [\mathcal{D},\mathcal{D}^{k-1}].
\]
As a module, $\mathcal{D}^k$ is generated by the set of vector
fields $V_I$, where $I$ describes the set of words with length
$k$. Moreover, due to Jacobi identity, we have
$[\mathcal{D}^i,\mathcal{D}^j] \subset \mathcal{D}^{i+j}$. This
flag is called the \textit{canonical flag} associated with the
differential system $\mathcal{D}$. H\"ormander's strong condition,
which we supposed to hold, states that for each $x \in
\mathbb{R}^n$, there is a smallest integer $r(x)$ such that
$\mathcal{D}^{r(x)}=\mathbb{R}^n$. For each $x \in \mathbb{R}^n$,
the canonical flag induces a flag of vector subspaces,
\[
\mathcal{D} (x) \subset \mathcal{D}^2(x) \subset \cdots \subset
\mathcal{D}^{r(x)} (x) =\mathbb{R}^n.
\]
The integer list $\left( \dim \mathcal{D}^k (x) \right)_{1 \leq k
\leq r(x)}$ is called the growth vector\index{Growth vector} of
$\mathcal{V}$ at $x$. The point $x$ is said to be a regular
point\index{Regular point} of $\mathcal{V}$ if the growth vector
is constant in a neighbourhood of $x$. Otherwise, we say that $x$
is a singular point. On a Carnot group, due to the homogeneity,
all points are regular.

Let $\mathcal{V}_i =\mathcal{D}^i / \mathcal{D}^{i-1}$ denote the
quotient differential systems, and define
\[
\mathcal{N} (\mathcal{D})=\mathcal{V}_1 \oplus \cdots \oplus
\mathcal{V}_k \oplus \cdots.
\]
The Lie bracket of vector fields induces a bilinear map on
$\mathcal{N} (\mathcal{D})$ which respects the grading:
$[\mathcal{V}_i,\mathcal{V}_j] \subset \mathcal{V}_{i+j}$.
Actually, $\mathcal{N} (\mathcal{D})$ inherits the structure of a
sheaf of Lie algebras. Moreover, if $x$ is a regular point of
$\mathcal{D}$, then this bracket induces a $r(x)$-step nilpotent
graded Lie algebra structure on $\mathcal{N} (\mathcal{D})(x)$.
Observe that the dimension of $\mathcal{N} (\mathcal{D})(x)$ is
equal to $n$ and that from the definition, $(V_1 (x),\ldots,V_d(x))$
Lie generates $\mathcal{N} (\mathcal{D})(x)$.
\begin{definition}
If $x$ is a regular point of $\mathcal{D}$, the $r(x)$-step
nilpotent graded Lie algebra $\mathcal{N} (\mathcal{D})(x)$ is
called the nilpotentisation of $\mathcal{D}$ at $x$. This Lie
algebra is the Lie algebra of a unique Carnot group which shall be
denoted $\mathbf{Gr}(\mathcal{D})(x)$ and called the tangent space
to $\mathcal{D}$ at $x$. The integer $D=\sum_{k=1}^{r(x)} k \dim
\mathcal{D}^k(x)$ is called the homogeneous dimension of
$\mathbf{Gr}(\mathcal{D})(x)$.
\end{definition}

\begin{proposition}
\label{small time density} Let $x$ be a regular point of
$\mathcal{D}$. Let $p_t$, $t>0$, denote the density with respect
to the Lebesgue measure of the solution of the stochastic
differential equation
\begin{equation}\label{SDEhypoelliptic}
X^x_t=x+\sum_{i=1}^d \int_0^t V_i (X_s^x)dB^i_s,
\end{equation}
where $(B_t)_{t \ge 0}$ is a $d$-dimensional fractional Brownian
motion with Hurst parameter $H>1/2$. We have,
\[
p_t (x) \sim_{t \rightarrow 0} \frac{C(x)}{t^{HD(x)}},
\]
where $C(x)$ is a non negative constant and $D(x)$ the homogeneous
dimension of the tangent space $\mathbf{Gr}(\mathcal{D})(x)$.
\end{proposition}

\begin{proof}
Let us first introduce some notations: For $k \ge 1$, we denote by 
$\Delta^k [0,t]$ the simplex of ordered $k$-tuples with values in $[0,t]$, \ie
\item
\[
\Delta^k [0,t]=\{ (t_1,\ldots,t_k) \in [0,t]^k, t_1 < \ldots < t_k \};
\]
If $I=(i_1,\ldots, i_k) \in \{1,\ldots,d\}^k$ is a word with length
$k$, we define the corresponding iterated integral of $B$ by
\begin{equation*}
\int_{\Delta^k [0,t]}  dB^I= \int_{0 < t_1 < \ldots < t_k < t}
dB^{i_1}_{t_1}  \cdots dB^{i_k}_{t_k}\;,
\end{equation*}
where the right hand side consists of nested Riemann-Stieltjes integrals.

We denote $\mathfrak{S}_k$ the group of the permutations of the
index set $\{1,\ldots,k\}$ and if $\sigma \in \mathfrak{S}_k$ and $I$ is a word $I=(i_1,\ldots,i_k)$,  
we denote by $\sigma \cdot I$ the word
$(i_{\sigma(1)},\ldots,i_{\sigma(k)})$.

If $\sigma \in \mathfrak{S}_k$, we denote $e(\sigma)$ the
cardinality of the set
\[
\{ j \in \{1,\ldots,k-1 \} , \sigma (j) > \sigma(j+1) \}\;,
\]
\ie $e(\sigma)$ is the number of raising sequences of $\sigma$.
Finally, if $I=(i_1,\ldots,i_k) \in \{ 1,\ldots, d \}^k$ is a word
\[
\Lambda_I (B)_t=\sum_{\sigma \in \mathfrak{S}_k} \frac{\left(
-1\right) ^{e(\sigma )}}{k^{2}\left(
\begin{array}{l}
k-1 \\
e(\sigma )
\end{array}
\right) } \int_{\Delta^k [0,t]} dB^{\sigma^{-1} \cdot I}.
\]

As a consequence of Proposition 23 in \cite{FriVic06Euler}, we get
the following approximation result
\[
X_t^{x}=\left[ \exp \left( \sum_{k = 1}^{r(x)}
\sum_{I=(i_1,\ldots,i_k)} \Lambda_I (B)_t V_I \right) \right](x)+
t^{\frac{r(x)+1}{2}} \mathbf{R} (t), \text{ } t \geq 0\;.
\]
Here, the remainder term $\mathbf{R}(t)$ is such that there exist $\alpha, c>0$ such that, for all $A
>c$,
\[
\lim \sup_{t \rightarrow 0} \mathbb{P} \left( \sup_{0 \leq s \leq
t} s^{H(N+1)} \mid \mathbf{R}_N (s) \mid \geq A t^{H(N+1)} \right)
\leq \exp \left( - \frac{A^{\alpha}}{c} \right).
\]
It is shown in \cite[Proposition 3.7, pp 66]{Bau2005} that
one can write
\[
\left[ \exp \left( \sum_{k = 1}^{r(x)} \sum_{I=(i_1,\ldots,i_k)}
\Lambda_I (B)_t V_I \right) \right](x)=F(B_t^{\ast} ),
\]
where $(B^{\ast}_t)_{t \geq 0}$ is the lift of $(B_t)_{t \geq 0}$
in the free Carnot group $\mathbb{G}_{r(x)} (\mathbb{R}^d)$ and
where $F$ is a map
\[
\mathbb{G}_{r(x)} (\mathbb{R}^d) \rightarrow \mathbb{R}^n
\]
that is Pansu differentiable at $1_{\mathbb{G}_{r(x)}}$. Furthermore, its
Pansu's derivative
\[
d_P F (1_{\mathbb{G}_{r(x)}}):\mathbb{G}_{r(x)}(\mathbb{R}^d)
\rightarrow \mathbf{Gr}(\mathcal{D})(x)
\]
is a surjective Carnot group morphism. Now, we have
\[
\left( F(B_t^{\ast} ) \right)_{t \geq 0} \eqlaw \left(
F\left(\Delta_{\sqrt{t}}^{\mathbb{G}_{r(x)}(\mathbb{R}^d)}
B_1^{\ast} \right)\right)_{t \geq 0},
\]
and
\[
F\left(\Delta_{\sqrt{t}}^{\mathbb{G}_{r(x)}(\mathbb{R}^d)}
B_1^{\ast} \right) \approx x
\Delta_{\sqrt{t}}^{\mathbf{Gr}(\mathcal{D})(x)} d_P F
(0_{\mathbb{G}_{r(x)}(\mathbb{R}^d)}) (B_1^{\ast}),
\]
in small times. Since $d_P F
(0_{\mathbb{G}_{r(x)}(\mathbb{R}^d)})$ is surjective, the random
variable
\[
d_P F (0_{\mathbb{G}_{r(x)}(\mathbb{R}^d)}) (B_1^{\ast})
\]
admits a density with respect to the Lebesgue measure on $\R^n$. This
concludes the proof.
\end{proof}

\bibliographystyle{Martin}
\markboth{\sc \refname}{\sc \refname}
\bibliography{referenceshormander}

\end{document}